\documentstyle[11pt]{article}
\begin{document}
\baselineskip=13pt plus 2pt minus 1pt
\bibliographystyle{plain}
\normalsize
\title{Exact sequences of semistable vector
bundles on algebraic curves}
\author{E. Ballico - B. Russo}
\def\length{\mbox{ length }}
\def\Supp{\mbox{ Supp}}
\def\End{\mbox{ End}}
\def\Ker{\mbox{ Ker}}

\def\parn{\par\noindent}
\newcommand{\Lm}{\mbox{Im}}
\def\rk{\mbox{rk}}
\def\Hom{\mbox{Hom}}
\newtheorem{thm}{Theorem}[section]
\newtheorem{prop}[thm]{Proposition}
\newtheorem{definition}[thm]{Definition}
\newtheorem{corollary}[thm]{Corollary}
\newtheorem{lemma}[thm]{Lemma}
\newtheorem{remark}[thm]{Remark}
\def\coker{\mbox { coker }}\def\N { \mbox{I}\!\! \!
\mbox{ {\bf N}}}
\def\Z{\bf Z}

\maketitle

\begin{abstract}
\noindent
 Let $X$ be a smooth complex projective curve of
genus $g\geq 1$. If $g\ge2$, assume further that $X
$ is either bielliptic or with general moduli. Fix
integers $r,$ $s,$ $a,$ $b$ with $r>1,$ $s>1$ and
$as \le br.$ Here we prove the existence of an exact
sequence

$$                              0 \to H \to  E
\to   Q \to 0
$$ of semistable vector bundles on  $X$ with
$\rk(H) = r,$ $\rk(Q) = s,$ $\deg(H) = a$ and
$\deg(Q) = b.$
\end{abstract}

\footnotetext{\parn 1991  Mathematics Subject
Classification  14H60.\parn This research was
partially supported by MURST and GNSAGA of CNR
(Italy). Both authors are members of the VBAC
Research group of Europroj. }
\normalsize\section*{\normalsize{      
\hspace{50mm}                      {\bf 
Introduction}}} Let $X$ be a smooth projective
curve on an algebraically closed field $k$. A vector bundle $E$ is said stable (resp. semistable) if for every proper subbundle $A$ of $F$ it holds $\mu(A)<\mu(E)$ (resp. $\mu(A)<\mu(E)).$ If MCD$(\rk(E),\deg(E))=1$ then $E$ is stable if and only if it is semistable. Since almost forty years these bundles are an active subject of study. In this paper we present a method to construct  a new semistable bundle, $E,$ as representative of an extension of the type
\begin{equation}
                0 \to H \to   E \to    Q \to  0
\label{eqm}\end{equation}
with $H$ and $Q$ semistable.
As one can easily see, a  necessary
condition for the stability (resp. semistability)
of
$E$ fitting in the exact sequence (\ref{eqm}) is
 $\mu(H) < \mu(Q)$ (resp. $\mu(H)   \leq
\mu(Q)). $ If $\rk(H) = \rk(Q)
 = 1$ it is known that these are also sufficient
conditions (see e.g. the properties of rank 2
stable vector bundles given in \cite{7}). Our main
results are the following theorems.

\begin{thm} {Let  $X$ be an elliptic curve. Fix
integers
$r,$ $s,$
$a,$ $b$ such that $r
\geq 1,$  $s\geq 1$ and $as  \leq br.  $ Then
there exists an exact sequence (\ref{eqm}) of
vector bundles on  $X$ with $\ rk(H) = r,$ $\rk(Q)
= s,$ $\deg(H) = a,$ $\deg(Q) = b,$ and $H,$  $E$
and  
$Q$ semistable.}\label{t01}
\end{thm}

\begin{thm} {Let  $X$ be a smooth complex
projective curve of genus $g\geq 2$ such that there
is a degree $2$ map $f:X  \to Y$ with $Y$ an
elliptic curve. Fix integers $r,$ $s,$ $a,$ $b$
such that $r
\geq 1,$  $s\geq 1$ and $as  \leq br.  $ Then
there exists an exact sequence (\ref{eqm}) of
vector bundles on  $X$ with $\ rk(H) = r,$ $\rk(Q)
= s,$ $\deg(H) = a,$ $\deg(Q) = b,$ and $H,$  $E$
and  
$Q$ semistable.}\label{t02}
\end{thm}

\begin{thm} {Let  $X$ be a general smooth complex
projective curve of genus $g\geq 2.$
 Fix integers $r,$ $s,$ $a,$ $b$ such
that $r\geq 1,$ $s\geq 1$ and
$as \leq br.$ Then there exists an exact sequence (\ref{eqm}) of
vector bundles on  $X$ with
$\rk(H) = r,$ $\rk(Q) = s,$ $\deg(H) = a,$ $\deg(Q)
= b,$ and $H,$  $E$ and   $Q
$ semistable.}\label{t03}\end{thm}

In section $1$ we will give a few general lemmata
and deduce Theorem \ref{t03} from Theorem
\ref{t02}. In section 2 we will prove Theorem
\ref{t01}. In section 3 we will prove
Theorem \ref{t01}.
\parn
{\em Acknowledgment}: we would like to thank professor P.E. Newstead for his
important suggestions and his kind support he gave us for all the period we
wrote this paper. The second author would like to thank the University of
Liverpool where part of the work was done.

\par
\vspace{3mm}
\section {\normalsize  \bf  {General remarks}} We
work in characteristic $0;$ in this case we may use the fact that, if $h:U \to
W$ is a finite morphism of smooth curves and   $F$ is a
semistable vector bundle on $W$, then $h^*(F)$ is
semistable.
\begin{remark}  {\rm Let  $X$ be a smooth
projective curve of genus $g$. It is well known
that, for any bounded family of vector bundles on
$X$ of fixed rank and degree, there exists a family
of vector bundles on
$X$ containing all members of the given family
and parametrised by an irreducible variety. For reader's sake we will give the
following short proof; up to a twist by a line bundle we may assume that every
such bundle is spanned; since $\dim(X) = 1,$ any such $r-$bundle, $E,$ is
spanned by
$r+1$ global sections, i.e. there is a surjection $h: {\cal O}_X ^{r+1} \to E;$ the
vector bundle
$\Ker(h)$ is isomorphic to $\det(E)^*$ and hence any such bundle $E$ appears as
cokernel of a family of embeddings of a family of line bundles of fixed degree into the trivial
bundle of rank
$r+1 .$  If
$g>0$ there exist semistable vector bundles on $X$
of any rank and any determinant. Hence, by the
openness of semistability, it follows that every
vector bundle on 
$X$ is the flat limit of a family of semistable
bundles with the same
determinant. If $g\ge2$, the
same holds with `semistable'
replaced by `stable'.\label{rr}}\end{remark}
\begin{lemma} Fix a flat family $\{X(t),\,\,\, t\in
T\}$ of smooth genus $g$ curves and flat families
$H(t),$ $Q(t)$ of vector bundles on the family
$\{X(t),
\,\,\,t\in T\}$ with $T$ integral. Assume the
existence of $0\in T$ such that ${\rm
Hom}(Q(0),H(0))=0
$. Then any extension of $Q(0)$ by $H(0)$ is the
flat limit of extensions of $Q(t)$ by $H(t)$ for $t$
in some open subset of $T$ containing $0$.
\label{lii}\end{lemma}
\parn
{\it Proof.} By semicontinuity we can find an open
subset $T_0$ of $T$ containing $0$ such that
$\Hom(Q(t),H(t))=0$ for all $t\in T_0$. Hence by
Riemann - Roch $h^1(Hom(Q(t),H(t)))$ is independent
of $t\in T_0$. It follows that there exists a
vector bundle $V$ over $T_0$ whose fibre over $t$
can be identified with $H^1(Hom(Q(t),H(t)))$;
moreover $V$ is the base of a family parametrising
all extensions of 
$Q(t)$ by $H(t)$ for $t\in T_0$ (see \cite{6}, Cor. 3.3 ). The result follows.
\begin{corollary}{Fix a flat family of curves as in
Lemma \ref{lii} and let
$H$,
$Q$ be vector bundles on
$X$ such that
${\rm Hom}(Q,H)=0$. Then any extension (\ref{eqm}) 
is the limit of a flat family of extensions
\[0\to H(t)\to E(t)\to Q(t)\to0\]with $H(t)$ and
$Q(t)$ semistable. Moreover, if $E$ is
semistable, we can suppose that $E(t)$ is
semistable.}
\label{c13}\end{corollary}
{\it Proof.} This follows from remark \ref{rr},
Lemma \ref{lii} and the openness of semistability. 
\parn
{\it Proof of Theorem \ref{t03} assuming Theorem
\ref{t02}.} If $as=br$, we can choose $H$ and $Q$
semistable and then take $E=H\oplus Q$. If $as<br$,
choose a flat family of curves as in Lemma
\ref{lii} such that
$X(t)$ is a general curve of genus $g$ for $t\ne0$,
while
$X(0)$ is bielliptic. Choose an extension
(\ref{eqm}) on $X(0)$ and note that, since $H$ and
$Q$ are semistable and $\mu(H)<\mu(Q)$, we have
$\Hom(Q,H)=0$. The result now follows by corollary
\ref{c13}.
\begin{lemma} Let $F$ and $G$ be vector bundles on
$X$ and set  $$\begin{array}{ll}t:=\max _{f\in
{\rm Hom}(F,G)}{\rm rk} (\Lm (f))  & h:= \max
_{f\in {\rm Hom}(F,G), {\rm rk} (Im(f))=t} \deg
(\Lm (f)).
\cr \end{array}$$ Then the subset of ${\rm
Hom}(F,G)$
$$ U:= \{ f\in {\rm Hom}(F,G)     \, \, |
\, \, {\rm rk} (\Lm (f))
=t, \deg (\Lm (f))=h\}$$ is open and dense
in ${\rm Hom}(F,G)$ and $\det(Im(f))$ is constant
for $f\in U$.
\label{lp}
\end{lemma}
{\it Proof.} Let  $ U':= \{ f\in
{\rm Hom}(F,G)   \, \, | \, \, \rk (Im (f)) =t
\}$. By the semicontinuity of the function  `rank', $U'$ is open and
dense in ${\rm Hom}(F,G)$. Now, the restriction to $U'$ of the function `degree' is semicontinuos;
hence $U$ is open and dense in ${\rm
Hom}(F,G)$. Since
${\rm Hom}(F,G)$  is rational,
 the morphism from $U$ to Pic$^h(X)$ mapping $f$
into $\det(Im(f))$ is constant.

\par
\vspace{3mm}
\section {\normalsize  \bf  {Elliptic curves;
proof of Theorem \ref{t01}}}
We suppose throughout this section
that $X$ is an elliptic curve. For
general facts about vector bundles on
$X$, see \cite{1}.\begin{remark} {\rm Let
$E$ be a vector bundle on the elliptic curve $X$
with
$\deg(E)>0$. Then, by Riemann-Roch,
$h^0(E)>0$.}\label{r21}\end{remark}
\begin{remark}{\rm Recall that
a polystable vector bundle on a smooth curve is a
direct sum of stable vector bundles with the same
slope. Note that, for any
rank $r$ and degree $d$,
there exist polystable
bundles on $X$. In fact, if
$(r,d)=l$, it follows from
the results of \cite{1}, theorem 7 and 10, pg. 433 and 442,  that there exist non-isomorphic stable bundles
$E_1,\ldots,E_l$ of rank $r\over l$ and degree
$d\over l$. Then
$E_1\oplus\ldots\oplus E_l$
is polystable and no two among its indecomposable factors are isomorphic. . Moreover
$E_1,\ldots,E_l$ can be
chosen  such that 
$E_1\oplus\ldots\oplus E_l$
has any given determinant. Note
also that any subbundle of
$E$ with the same slope as
$E$ is a  direct factor of it.}\label{r22}\end{remark} In
order to prove theorem
\ref{t01} we need the following
fundamental lemma on polystable vector bundles.

\begin{prop}{Let $F$, $G$ be polystable vector
bundle on an elliptic curve $X$
 with $\rk(F) \geq \rk(G)$ and $\mu(F) < \mu(G).$
Assume that no two among the indecomposable factors
of   $F$ (resp. of $G)$ are isomorphic. Let $f\in
U$, where $U$ is defined as in Lemma \ref{lp}. Then
\par
        (a) We have $\rk(\Lm(f)) = \rk(G).$
\par
        (b) If $\rk(F) > \rk(G),$ then $f$ is
surjective.}
\label{pe}

\end{prop}Let us notice that (a) and (b) have the
following `dual' version
\par (a$'$) If $\rk(F)\leq \rk(G)$ a general $f\in
{\rm Hom}(F,G)$ is injective.
\par (b$'$) If $\rk(F) < \rk(G)$ then a general 
$f\in {\rm Hom}(F,G)$ is injective and $\coker
(f)$ is torsion free.

\vspace{3mm}\parn {\it Proof.}  We will use the same
notation as in Lemma \ref{lp}. Since $Hom(F,G)$ has
slope
$\mu(G) - \mu(F) > 0,$ we have
$\Hom(F,G)
 \ne 0$ by remark \ref{r21}. 
\par First we prove that if (a) is
 true  then (b) follows. Let us fix $F$ and an
integer $s.$ Then consider vector
 bundles $G$ of rank $s$. Since (a) holds we may
suppose $\rk(\Lm(F))=\rk(G)=s$. Then choose $d
_0\in\Z$ such that $(d _0-1)
/s \le
\mu(F)<d_0/s$. We prove first that (b) holds if
$\deg(G)=d_0$. In this case, we have, since $F$ is
polystable,
$$ {(d _0-1)\over
\rk(G)}\le\mu(F)\leq
\mu(\Lm(f)) \le {d _0\over \rk (G)}.$$ If
$\mu(\Lm(F)=d_0/\rk(G)$, we are done.
Otherwise, we have  $(d_0-1)/\rk(G) =
\mu(F)=\mu(\Lm(F))$; hence by remark \ref{r22}
$\Lm(f)$ is a
 direct factor of $F$, say   $F   \simeq 
\Lm(f)\oplus B.$  By the
assumption on the indecomposable factors of   $F$
we have $\Hom(B,\Lm(f))=0$, while by remark
\ref{r21} $\Hom(B,G)\ne0$.
Hence there is $u: B   \to G$ with $\Lm(u)$ not
contained in $\Lm(f).$ Define $v:F \to G$ with $v|B
= u$ and sending the factor $\Lm(f)$ of   $F$ to
the subsheaf $\Lm(f)$ of $G$ via the identity. By
construction $\Lm(v)$ strictly contains  $\Lm(f).$
Since $\rk (G)= s =\rk (\Lm (v))$ then $\deg (\Lm
(v) )> \deg (\Lm (f)) $ which is
 absurd because $ f \in U$. This proves (b) in
the case $\deg(G)=d_0$.
To complete the proof of (b), we check by
induction that the set
$$\begin{array}{l}B:=\{d\in \N \,\,\,| \,\,\,d\geq
d _0 \mbox{ and }\forall \,\, \mbox{ polystable
vector bundle } G \mbox{ of rank } s \cr\mbox{ and
degree } d,
\mbox{ a  general }
  f: F \to G \,\,\, \mbox{ is
surjective}\}.\end{array}$$ is the set of all the
integers greater than or equal to $ d_0$. By the
previous argument $d _0\in B$. Then we assume
$d-1\in B$ and prove that
$d\in B.$ Suppose $d\not\in B$. By the inductive
step, for any
$L
\in
\mbox{Pic}^0(X)$ and for every polystable vector
bundle, $T,$ of rank $s$ and degree $d-1$  the
general map $F\to T\otimes L$ is surjective. Since
(a) and (a') hold then  the general map $T\otimes L
\to G$ is injective. Then the composition $f
_L:F\to G$ has image of rank $s$ and degree $d-1.$
If $h$ is the maximum defined in lemma \ref{lp},
then $h = d-1$ and hence  $f _L\in U$. Now we can
vary $L\in \mbox{Pic}^0(X)$ (which is infinite)
and  find images of maps $f _L\in U$ with
non-isomorphic determinants.  This contradicts
Lemma \ref{lp}.\par Finally let us prove part (a).
 We will use induction on  $\rk(G).$ If $\rk(G) =
1,$ the first assertion is obvious. Suppose
$\rk(G)>1$ and the inductive step. In order to
obtain a contradiction we assume $t < \rk(G).$ Now
the proof is divided in two cases. First assume
${h\over t} > \mu(F)$.  By remark \ref{r22} we can
fix a polystable bundle  $M$ with $\rk(M) = t,$
$\deg(M) = h$ and factors which are mutually
non-isomorphic. Then  by the inductive step, for any
$L\in
\mbox{Pic}^0(X)$ the rank of the image of the
general map  $F\to M\otimes L$ is $t.$ Since (a)
implies (b), the image is, in fact, all of
$M\otimes L.$ If ${h\over t} < \mu(G)$, using  the
dual form of the statement of the theorem, the
inductive step  applies to the pair $(M
 \otimes L ,G).$ Then we get by composition a map
$f _L :F \to G$ which has image in $G$
with maximal degree and rank.  Varying $L\in
\mbox{Pic}^0 (X)$ we get a contradiction. If
${h\over t} = \mu(G)$, we can choose $f\in U$
so that $\Lm(f)$ is a direct factor of $G$. Let
$l=(h,t)$. Then we can write
$\Lm(f)=G_1\oplus\ldots\oplus G_l$ and
$G=G_1\oplus\ldots\oplus G_m$ with $l<m$ and all
$G_j$ stable and non-isomorphic. By [\cite{1}, theorem 7 and 10, pg. 434 and 442], the bundles $\det(G_j)$ are all
non-isomorphic
. Hence $G'=G_1\oplus\ldots\oplus
G_{l-1}\oplus G_{l+1}$ has a different determinant
from $G_1\oplus\ldots\oplus G_l$. But, by the
inductive hypothesis there exists a surjective
homomorphism $f':F\to G'$. Since $f$ and $f'$ both
belong to $U$, this contradicts Lemma \ref{lp}.
Now assume 
${h\over t}=
\mu(F).$ Since   $F$ is polystable, then
$\Lm(f)$ is a direct factor of $F,$ say   $F  
\simeq  \Lm(f)\oplus B.$ By the assumption on the
indecomposable factors of   $F$ and the
semistability of its factors, we have $
\Hom(B,\Lm(f)) = 0,$ while $\Hom(B,G) \not=
 0$ by remark \ref{r21}. Hence there is $f': F \to
G$ with
$f'|B
\not= 0$ and $f'|\Lm(f) = f.$ Thus
 $\rk(f') > \rk(f)$, a
contradiction.\begin{corollary}  Let $F$ and $G$ be   polystable vector bundles with non-isomorphic factors, $\rk (F)=\rk (G)$ and $\deg (F)=\deg (G)+1$  Then $G$ is a subsheaf of $F.$ Furthermore a
general positive elementary transformation of  $G$
is semistable.
\label{r24}\end{corollary}\parn {\it Proof.} Apply the
previous proposition to the pair
 $(G,F).$ The last assertion follows from the
openness of semistability.
\begin{remark}{\rm
Applying twice corollary \ref{r24} one may obtain for every pair $G'$ and $F'$ of polystable vector bundles with non-isorphic factors, the  same rank and $\deg (G')-\deg(F')=2,$ that the generic map, $G'\to F',$ is an inclusion.}\label{r25}
\end{remark}
\parn {\it Proof of Theorem 0.1.} If $as=br$, take
any
$H$ and   $Q$ semistable and set $E := H\oplus Q.$
If $as<br$, by Proposition \ref{pe} there
exist polystable vector bundles $E',$
$Q'$ with non-isomorphic  factors and a surjection
$u: E'   \to Q'$ with $\rk(E') = r+s,$ $\deg(E') =
a+b,$
 $\rk(Q') = r,$ $\deg(Q') = b.$ To prove the
theorem it is sufficient to show that for general
$(E',Q',u)$, $\Ker(u)$ is semistable. Since $E'$
is semistable and $\mu(E') <
\mu(Q')$, $\Hom(Q',\Ker(u))=0$. The result follows
from corollary \ref{c13}.
\begin{remark}
{\rm Since the condition for a bundle to be polystable with non- isomorphic factors is preserved for general deformation, what we get at the end of the proof of Theorem 0.1 is an exact sequence in which the last two terms are polystable vector bundles with non isomorphic factors while the first one is semistable.}\label{intell}
\end{remark}

\section{ \normalsize {\bf Proof of theorem \ref{t02} } }

Let $X$ be a curve of genus $g$ and $f:X\to Y$ a $2:1$ morphism on an elliptic curve $Y.$ Let
$\sigma:  X   \to  X$ be the involution corresponding to the morphism $f,$ i.e. with
$Y  = X/\sigma.$ 
First let us notice that if $as = br,$ we can choose $H$ and $Q$ semistable and set $E:=
H\oplus Q.$ Hence from now on, we may assume $as < br.$
In the proof the
following remark will turn out to be quite useful.
\begin{remark}{\rm
Let $U$ be a vector bundle on $X.$ By the descent theory $U$ is of the form
$f^*(F)$ with $F$ a vector bundle on $Y$ if and only if the following two conditions are
satisfied:\par i) $U$ is $\sigma$ invariant,
\par
ii) for every ramification point $Q$ of $f$ $\sigma$ acts as the identity on the
 fiber $U_Q.$\par
Notice that a saturated subbundle of $f^*(F)$ satisfies condition ii). Hence a saturated subbundle
$A$ of $f^*(F) $ is of the form $f^*(B)$ with $B$ a vector bundle on $Y$ if and only if $A$ is
$\sigma$-invariant. Furthermore, the saturation in $f^*(B)$ of a $\sigma$-invariant subsheaf of
$f^*(B)$ has even degree.
 \label{ram}}\end{remark}
In order to prove the theorem we need the following fundamental lemmata.

\begin{lemma}
Let $M$ and $N$ be semistable vector bundles on $Y$ with $\rk (M)=\rk (N),$ $\deg (M)=\deg(N)+1$ and  $N\hookrightarrow M.$  Set $G=f^*(N)$ and $T=f^*(M).$ Then every vector bundle $E$ with $G\hookrightarrow E\hookrightarrow T$ and 
$\length (E/G) = 1 ,$  is semistable. 
\label{pl1}
\end{lemma}
\begin{lemma}
Let $M$ and $N$ be semistable vector bundles on $Y$ with $\rk (M)=\rk (N),$ $\deg (M)=\deg(N)+2$ and  $N\hookrightarrow M.$  Set $G=f^*(N)$ and $T=f^*(M).$ Then every vector bundle $E$ with $G\hookrightarrow E\hookrightarrow T$ and 
$\length (E/G) = 2 ,$  is semistable. 
\label{pl2}
\end{lemma}
The proof of lemma \ref{pl1} turns out to be a special case of the proof of lemma
\ref{pl2}. Several situations of the latter do simplify in the former (for example
$A$ is always saturated in the proof of lemma \ref{pl1}). Hence we will present only
the proof of lemma
\ref{pl2}.\parn
{\it Proof of \ref{pl2}.} Since the field characteristic is $0,$ then $G$ and $T$ are semistable. Let $E$ be a subsheaf of $T$ containing $G$ and with 
length$(T/E)=2.$ In order to find  a contradiction we assume that $E$ is not semistable. Then there
exists a stable proper subbundle $A $ of $E$
with
$\mu(A)>\mu(E).$
 Since $\mu(G) =
\mu(E) -2/(r+s)$ and $E''$ is semistable, we see that $A\cap E''  \not= A.$ Hence
there exists an inclusion $A\hookrightarrow T.$ 
Moreover note that $A  \not= 
\sigma(A)$ otherwise $A=f^*(N)$ would be contained in $E\cap \sigma(E)$ which is
impossible because the latter is contained in $E''.$
Now our intent is to prove that $A\cap \sigma(A)$ and $A+\sigma(A)$ are saturated
in $T.$ 
In order to prove this, we distinguish two cases depending if  $A$ is saturated in $T$
or not. 

If $A$ is saturated in $T$ one has $A\cap \sigma(A)$ saturated in $T.$ In fact the
 saturation $F$ of $A\cap \sigma(A)$ is contained in $A$ and since it is
$\sigma$-invariant, it is also contained in $\sigma(A).$ Then the saturation
coincides with $A\cap\sigma(A).$  Anyway in this case one might still have
that
$A+\sigma(A)$ is not saturated in
$T.$ Then one observes that by the saturation of $A\cap \sigma(A),$ $\deg
(A+\sigma(A))$ is even. Hence setting 
$K$ the saturation of $\deg (A+\sigma(A)),$
one gets

{\small
$$\begin{array}{l}
\mu(K)\rk (K)=\deg (K)\geq \deg
(A+\sigma(A))+2=\mu(A+\sigma(A))\rk(A+\sigma(A))  +2\\ \\
\geq\mu(A)\rk(A+\sigma(A))+2
>\mu(E)\rk(A+\sigma(A))+2=\mu(T)\rk(A+\sigma(A))-\\ \\
-2{\rk(A+\sigma(A))\over \rk (E)} +2,
\end{array}$$}
which contradicts the semistability of $T.$

Suppose now that $A$ is not saturated in $T.$ Let $A''$ be its saturation in $T.$
Call $A'$ the first term of the Harder-Narasimhan filtration of
$A''.$ If $A''$ is not semistable, then $A'$ is a proper subbundle of $A''.$
 In any case we get a semistable subbundle
$A'$ of
$T$ of slope
$\mu(A')\geq\mu(A)+{1\over \rk (A)}>\mu(E)+{1\over
\rk (A)}.$  Then  the inequalities

{\small $$\begin{array}{l}\mu(K)\rk (K)=\deg (K)\geq \deg(A'+\sigma(A'))+1\geq
\mu(A)\rk(A+\sigma(A)) +\\ \\+ {\rk(A+\sigma(A))\over \rk (A)} +1
>\mu(T)\rk(A+\sigma(A))-2{\rk(A+\sigma(A))\over
\rk (E)}+{\rk(A+\sigma(A))\over \rk (A)}+1\\ \\\geq   \mu(T) \rk (A+\sigma(A)) 
+{2\over
\rk (E)}(\rk (E)-\rk (A))\end{array}$$}
contradict the semistability of $T$ and they imply the saturation
of $A+\sigma(A).$ Hence degree of $A\cap \sigma(A)$ is even. Therefore  

{\small $$\begin{array}{l}\mu(F)\rk (F)= \deg(F)\geq
\deg(A'\cap\sigma(A'))+2=2\mu(A')\rk (A)-\\ \\ -\mu(A'+\sigma(A'))\rk(A+\sigma(A))+2
\geq 2\mu(A)\rk(A)+2 -\mu(T)\rk(A+\sigma(A)) +\\ \\+2>
2\mu(E)\rk (A)-\mu(T)\rk(A+\sigma(A)) +4=\mu(T)\rk(A\cap\sigma(A))-{4\rk (A)\over
\rk (E)}+4 \end{array}
$$}
which again contradicts the semistability of $T$ implying the saturation of $A\cap
\sigma(A).$
At this point we know that even if 
 $A$ is not saturated we  have
$A'\cap \sigma(A')\simeq f^*(U),$ $A'+ \sigma(A')\simeq f^*(B).$ 
Observe that $A+\sigma(A)$ is saturated whenever $A\cap \sigma(A)$ is zero or
not. Once we know this we can prove that $A\cap \sigma(A)\not =0,$ indeed.
In fact in order to prove a
contradiction we assume that 
$\rk(A+
\sigma(A)) = 2\rk(A).$  Thus $f^*(B)$ splits. We may assume $B$ indecomposable 
because
$\rk(A)$ is minimal. Note that the branch locus $B(f)$ of $f$ is not empty because $g > 1.$ By the
projection formula we have $f^*(\End(f^*(B)))    \simeq  \End(B)\oplus \End(B)(-B(f)).$ Since
$\deg(B(f))>0$ and End$(B)(-B(f))$ is semistable of degree less then zero, we have
$H^0(f^*(\mbox{End}(f^*(B)))\simeq H^0(\mbox{End}(B)),$ contradicting the indecomposability of $B.$
\vspace{2mm}\parn
{\bf  Claim:} We have $h^0(X,End(f^*(B/U))) = h^0(Y,End(B/U)).$\parn
Assume the claim 
and consider the exact sequence $$ 0\to A\cap \sigma(A)\to A+\sigma(A) \to A+\sigma(A)/A\cap
\sigma(A).$$ By construction $f^*(B/U)=A+\sigma(A)/A\cap
\sigma(A)$ splits into two direct factors $A/A
 \cap \sigma(A)$ and  $\sigma(A)/A \cap \sigma(A).$ 
 The projection of $f^*(B/U)$ into its factor $A/A \cap \sigma(A)$ does not come
 from an element of $H^0(Y,End(B/U)).$ Hence  we get a contradiction and the theorem
is proven.\parn
{\it Proof of the Claim:}  By the projection formula we have
$$h^0(X,End(f^*(B/U))) = h^0(Y,End(B/U)) + h^0(Y,End(B/U)(-B(f))).$$ 
Hence it is sufficient to show the
 vanishing of $h^0(Y,End(B/U)(-B(f))).$ Since $Y$ is elliptic and card$(B(f)) =
2g-2$ by Riemann - Hurwitz theorem, it is sufficient to show that either $B/U$ is indecomposable or
that the difference, $t,$ between the highest and the lowest slope of its direct factors is less
then $ 2g-2.$ Suppose $B/U$ decomposable. Let us prove that $t < 2$ and hence $t < 2g-2$ for $g
\geq 2.$ In order to obtain a contradiction we assume the existence of a direct factor $A'$ of $B/U$
of slope $\geq \mu(A)/2 + 2.$ Define $B'$ by $B/U = A'\oplus B'.$ Since $A$ is stable and
$f^*(B/U)$ is a quotient of $A\oplus \sigma(A),$ every direct factor of $f^*(B/U)$
has slope $\geq
\mu(A).$ Hence our assumption implies that
$\deg(B/ U) \geq
\mu(A)\rk(A')/2 + 2\rk(A') +
\mu(A)\rk(B')/2 \geq \mu(A)\rk(B/U)/2 + 2.$
Moreover 

{\small $$\begin{array}{l}\deg(A+ \sigma(A))=\deg(f^*(B/U))+\deg(A\cap \sigma(A)) =
\deg(f^*(B/U))+ \\ \\ +2\deg (A) - \deg(A+\sigma(A)).\end{array}$$} Hence
{\small $$\begin{array}{l}\deg(A+
\sigma(A))=\deg(f^*(B/U))/ 2+ \deg (A)\geq
\mu(A)\rk(B/U)/ 2+ \\ \\+ 2+ \mu(A)\rk(A)=\mu(A)(\rk(A+\sigma(A)) +2. \end{array}$$}
Since
$\mu(A)>
\mu(T) - 2/\rk (T),$ and $T$ is semistable, we obtain a contradiction of the semistability of $T.$ Hence every
direct factor $A''$ of $B/U$ has $\mu(A)/2   < \mu(A'') <
\mu(A)/2 + 2.$ Hence $t < 2,$ as
 wanted.
\vspace{2mm}\par

\begin{corollary}
Let $N$ be a polystable bundle on $Y.$ Set $G=f^*(N).$   Then there exists a
semistable bundle $E$  containing $ G$ with
$\length (E/G) = 1 .$ 
\label{prec1}
\end{corollary}
{\it Proof of \ref{prec1}.} By remark \ref{r22} and corollary \ref{r24}, we can find a polystable bundle $M$ on
$Y,$ with non-isomorphic factors, containing  $N, $ with length$(M/N)=1.$  Hence we
can find a semistable bundle $T=f^*(M)$ and an inclusion
$G\hookrightarrow T$ compatible with $\sigma$ and with length$(T/G)=2.$  Set
Supp$(T/G)  = \{P,\sigma (P) \}.$ Since semistability is an open
condition, we may assume $
\sigma (P)
\not = P.$ Then there exists a vector bundle $E$ with $G\hookrightarrow E\hookrightarrow T,$ Supp$(E/G)=\{P\}$ and Supp$(T/E)=\{\sigma(P)\}.$ Hence  we may apply lemma \ref{pl1}.
The same proof gives the following corollary
\begin{corollary}
Let $N$ and $M$ be  polystable bundles on $Y$ with $N\hookrightarrow M$ and length$(M/N)=1. $ Set $G=f^*(N)$ and $T=f^*(M).$   Then there exists a
semistable bundle $E$ with $G\hookrightarrow E\hookrightarrow T$ and 
$\length (E/G) = 1 .$ 
\label{prec2}
\end{corollary}
\begin{corollary}
Let $F'$ a polystable bundle on $Y.$ Set $T=f^*(F').$   Then there exists a
semistable bundle $E$  contained in $ T$ with
$\mbox{length} (T/E) = 1 .$   
\label{dual1}\end{corollary}
{\it Proof.} The statement is just the dual version of corollary \ref{prec1} or it
derives directly by lemma
\ref{pl1}.
\vspace{5mm}\parn
{\it Proof of Theorem \ref{t01}}: Unless otherwise specified, every exact sequence
of vector bundles on $X$ will consist of a subbundle of rank $r,$ a middle term of rank
$r+s$ and a quotient of rank $s.$

The proof is divided into four cases according to the parity of the two integers $a$ and
$b.$
\vspace{2mm}\parn
i) Here we assume $a$ even, $b$ even. By theorem 0.1, we can find an exact
sequence
\begin{equation}
                0   \to H'   \to E'   \to Q'   \to 0             \label{et01}
                \end{equation}
of semistable vector bundles on $Y $ with $\deg(H'
) = a/2,$ $\deg(Q') = b/2.$ Set $H:= f^*(H'),$ $E:= f^*(E')$ and $Q:= f^*(Q').$ The pull-back of (\ref{et01}) gives the desired exact sequence. 

\vspace{2mm}\parn
ii) Here we assume $a$ even, $b$ odd. By remark \ref{intell} and corollary \ref{r24} we have on $Y$
the following diagram 

$$\begin{array}{ccclcccccc}
&&&&&0&&0&&\\
&&&&&\downarrow&&\downarrow&&\\
0&\to& H'&&\to& F'&\to &Z'&\to& 0\\
&&\downarrow&\simeq &&\downarrow&&\downarrow&&\\
0&\to& H'&&\to& E'&\to &Q'&\to& 0\\
\end{array}$$
where $H'$ is a semistable vector bundle with $\deg (H')=a/2,$ while $F',$ $Z',$ $E',$ $Q'$ are
polystable bundles with non-isomorphic factors, $\deg (E')={a+b+1\over 2},$ $\deg (Q')= {b+1\over
2},$ $\deg (F')={a+b-1\over 2},$ $\deg (Z')= {b-1\over 2}.$ 
Pulling back this diagram and applying corollary \ref{prec2} to $j: F'\hookrightarrow E'$ and 
$i: Z'\hookrightarrow Q',$ one gets the desired exact sequence:

\begin{equation}
        0   \to H   \to E   \to Q   \to 0                   \label{ep1}
\end{equation}
of semistable bundles on  $X$ with  $\deg(H) = a,$
 $\deg(Q) = b.$ 
\vspace{2mm}\parn
iii) Here we assume $a$ odd, $b$ even. Just dualize the previous result.
\vspace{2mm}\parn
iv) Here we assume $a$ odd, $b$ odd.  By remark \ref{intell}, corollary \ref{r24},  remark \ref{r25} and the snake lemma for the injectivity of the first vertical map, we have on $Y$
the following diagram
$$\begin{array}{ccccccccc}
&&0&&0&&0&&\\
&&\downarrow&&\downarrow&&\downarrow&&
\\
0&\to&U'&\to &F'&\to& Z'&\to& 0\\
&&\downarrow&&\downarrow&&\downarrow&&
\\
0&\to& H'&\to &E'&\to&
Q'&\to& 0\\
\end{array}$$
where $H'$ and $U'$ are semistable bundles with  
$\deg (H')={a+1\over 2}$ and  $\deg (U')= {a-1\over 2},$ 
while $F',$ $Z',$ $E'$ and $Q'$ are polystable bundles with non-isomorphic
factors and degrees $\deg  (Z')={b-1\over 2},$ 
and $\deg (Q')={b+1\over 2}.$
Then pullback the diagram to $X.$ Use corollary \ref{prec2} for 
the vertical map on the right, finding a semistable bundle, $Q,$ with $f^*(Z')
\hookrightarrow Q\hookrightarrow f^*(z).$  

Now set $\Supp
(f^*(E')/f^*(F')) = \{P,Q,\sigma(P),\sigma(Q)\}.$ Then as in the proof of
corollary \ref{prec1}, there exists  a vector bundle $E$ with
$f^*(F')\hookrightarrow E \hookrightarrow f^*(E')$ with $\Supp(E/f^*(F')) =
\{P,Q\}$ and  $\Supp(f^*(E')/E)=\{\sigma(P),\sigma(Q)\}.$ Hence we may apply
lemma \ref{pl2}.

\vspace{5mm}\parn 
\def\barbara{\small Barbara Russo \parn Universit\`a degli studi di Trento \parn
 via sommarive, 14 \parn
38050 Povo (TN), Italy; \parn  e-mail: russo@alpha.science.unitn.it \parn fax: Italy+ 461-881624}
\def\edoardo{\small Edoardo Ballico \parn  Universit\`a degli studi di Trento \parn via sommarive,
14 \parn 38050 Povo (TN), Italy; \parn e-mail: ballico@alpha.science.unitn.it \parn fax:
Italy+ 461-881624}
\begin{tabular}{p{6.5cm}p{6.5cm}}
\edoardo & \barbara
\end{tabular}


\begin{thebibliography}{dd}


\bibitem[1]{1}  M. F. Atiyah,  `Vector bundles over an elliptic curve'  Proc. London Math. Soc. (3)
7 (1957), 414-452; reprinted in: Michael Atiyah Collected Works, Vol. 1, pp. 105-143.
\bibitem[2]{2}  C. Banica, M. Putinar and G. Schumaker,  `Variation der globalen
 Ext in Deformationen kompakter komplexer Raume',  Math. Ann.  250 (1980), 135-1
55.
\bibitem[3]{3}  A. Hirschowitz,  `Probl\'emes de Brill-Noether en rang sup\'erieur', unpublished
preprint partially printed as ref. [4] .
\bibitem[4]{4}  A. Hirschowitz,  `Probl\'emes de Brill-Noether en rang sup\'erieur', Comptes Rendus 
Acad. Sci. Paris  307 (1988), 153-156.
\bibitem[5]{5}  H. Lange,  `Zur Klassifikation von Regelmannigfaltigkeiten',  Math. Ann. 262
(1983), 447-459.
\bibitem[6]{6} H. Lange,  `Universal families of extensions', J. Algebra 83 (1983),
101-112.
\bibitem[7]{7} H. Lange and M. S. Narasimhan,  `Maximal subbundles of rank two
vector bundles on curves',  Math. Ann.  266 (1983), 55-72.
\end{thebibliography}
\end{document}